\documentclass{article}
\usepackage{amsmath, amssymb, amsthm}

\newtheorem {theorem}{Theorem}[section]

\newtheorem{corollary}[theorem]{Corollary}
\newtheorem{lemma}[theorem]{Lemma}
\newtheorem {proposition}[theorem]{Proposition}
\theoremstyle {definition}
\newtheorem {definition}[theorem]{Definition}
\theoremstyle {remark}
\newtheorem {remark}[theorem]{Remark}
\newtheorem {example}[theorem]{Example}
\newcommand{\fm}{\ensuremath{\mathfrak m}}
\newcommand{\fa}{\ensuremath{\mathfrak a}}
\newcommand{\un}{\ensuremath{\underline}}

\newcommand{\sop}{\text{system of parameters }}

\newcommand{\cm}{\text{Cohen-Macaulay }}
\newcommand{\dfil}{$M_0\subset M_1\subset \cdots \subset M_{t-1}\subset M_t=M$\ }

\def\Ann{\operatorname{Ann}}

\begin{document}
\title{dd-sequences and partial Euler-Poincar\'{e} characteristics of Koszul complex}
\author{Nguyen Tu Cuong\footnote{Email: ntcuong@math.ac.vn} \  and Doan Trung Cuong\footnote{Email: dtcuong@math.ac.vn}\\
Institute of Mathematics\\
18 Hoang Quoc Viet Road, 10307 Hanoi, Vietnam}
\date{March 23, 2005}
\maketitle
\begin{abstract} The aim of this paper is to introduce a new  notion of  sequences called  dd-sequences and show that this notion  may be convenient for  studying the polynomial property of  partial Euler-Poincar\'e characteristics of the Koszul complex  with respect to the powers of a system of parameters.  Some results about the dd-sequences, the partial Euler-Poincar\'e characteristics and the lengths of local cohomology modules are also presented in the paper. There are also applications of dd-sequences on the structure of  sequentially Cohen-Macaulay modules.
\end{abstract}

\section{Introduction}
\markright{Introduction}

Let $(R,\fm)$ be a commutative local Noetherian ring and $M$ a finitely generated $R$-module of dimension $d$. Let $\un x =(x_1,\ldots,x_d)$ be a \sop of $M$. Denote by $H_i(\un x;M)$  the $i$-th Koszul homology module of $M$ with respect to the \sop  $\un x$. Following Serre \cite [Appendice II]{serre},  the $k$-th Euler-Poincar\'{e} characteristic of $M$ with respect to   $\un x$ is defined by 
$$\chi_k(\un x;M)=\sum_{i=k}^\infty(-1)^{i-k}\ell(H_i(\un x;M))$$
for  $k=0, 1,\ldots , d$, where $\ell (N)$ stands for the length of the $R$-module $N$. Let $\un n =(n_1,\ldots,n_d)$ be  a $d$-tuple of positive integers. We denote by $\un x(\un n)$  the system of parameters   $(x_1^{n_1},\ldots,x_d^{n_d})$. Then we can consider $\chi_k(\un x(\un n);M)$  as a function in $\un n$.  It is well-known by Garcia Roig \cite {roig} that this function in general is not a polynomial. But V. T. Khoi and the first author have proved in \cite {ck} that the least degree of all polynomials in $\un n$ bounding above the function $\chi_k(\un x(\un n);M)$ is independent of the choice of the \sop $\un x$ and they denote this invariant by $p_k(M)$. 
It was shown in \cite {c3} (in that paper  $\chi_1(\un x(\un n);M)$ and $p_1(M)$ are denoted by $I_M(\un n ;\un x)$ and $p(M)$ respectively) that  the $1$-st Euler-Poincar\'{e} characteristic of $M$ is a polynomial in $\un n$ provided $\un x$ is a p-standard  \sop of $M$. Moreover,  in this case the polynomial has a very simple formula as follows 
$$\chi_1(\un x(\un n);M)=\sum_{i=0}^{p_1(M)}n_1\ldots n_ie(x_1,\ldots,x_i;(0:x_{i+1})_{M/(x_{i+2},\ldots,x_d)M}),$$
where, for a \sop  $\un y$ of a $R$-module $N$, we denote by $e(\un y; N)$ the Serre multiplicity of $N$ with respect to  $\un y$ and set $$e(x_1,\ldots,x_i;(0:x_{i+1})_{M/(x_{i+2},\ldots,x_d)M})=\ell ((x_2,\ldots , x_d)M:x_1/(x_2,\ldots , x_d)M)$$
 for $i=0$.
Recall that  a system of parameters $\un x=(x_1,\ldots,x_d)$ of $M$ is called a p-standard  \sop if $x_i\in \fa(M/(x_{i+1},\ldots,x_d)M)$ for $i=1,\ldots,d$, where $\fa(M)=\fa_0(M)\ldots\fa_{d-1}(M)$ and $\fa_i(M)$  is the annihilator ideal of the $i$-th local cohomology module $H_\fm^i(M)$ with respect to the maximal ideal $\fm$. It should be mentioned that the notion of p-standard  \sop  has been used as an important tool for the resolution of Macaulayfication by T. Kawasaki \cite {ta1, ta2}. 
Therefore it raises to a natural  question:  Is $\chi_k(\un x(\un n);M)$ still  a polynomial in $\un n$ for all $k>1$ provided  $\un x$ is a p-standard system of parameters?

In order to give an affirmative answer to the above question we introduce in this paper the concept of dd-sequence, which is a slight generalization of the notion  of p-standard sequence. Note here that while there  does not exist a  necessary and sufficient condition for a \sop to be p-standard, we can characterize a \sop to be a dd-sequence in term of the $1$-st Euler-Poincar\'{e} characteristic. Therefore we can solve the above question for  the class of systems of parameters which are dd-sequences. Moreover, as an application, we can study more details about  systems of parameters,  which are dd-sequences, as well as the partial Euler-Poincar\'{e} characteristics with respect to them.  

Let us to summarize the main results of this paper. The paper is divided into 6 sections.
We  need first  a necessary and sufficient condition for the $k$-th Euler-Poincar\'{e} characteristic $\chi_k(\un x(\un n);M)$ of $M$ to be a polynomial in $\un n$. Hence the following theorem, which is a generalization of a result of \cite {c1}, is the main result of Section 2.
\begin{theorem}\label{main3}
Let $M$ be a finitely generated $R$-module, $\un x=(x_1,\ldots,x_d)$ a system of parameters of $M$ and $n_0$  a positive integer. Then the $k$-th Euler-Poincar\'{e} characteristic  $\chi_k(\un x(\un n);M)$ is a polynomial in $\un n=(n_1,\ldots , n_d)$ for all $n_i\geq n_0$ if and only if  the following condition is satisfied 
$$(0:x_i^{n_i})_{H_{k-1}(x_1^{n_1},\ldots,x_{i-1}^{n_{i-1}},x_{i+1}^{n_{i+1}},\ldots,x_d^{n_d};M)} =(0:x_i^{n_0})_{H_{k-1}
(x_1^{n_1},\ldots,x_{i-1}^{n_{i-1}},x_{i+1}^{n_{i+1}},\ldots,x_d^{n_d};M)}$$
for all $i=1,\ldots,d$ and 
$n_1,\ldots, n_d\geq n_0 $.
\end{theorem}

Next, in  Section 3
we use the concept of d-sequence of Huneke (see \cite {ch}) to  introduce  the notion of dd-sequences. A sequence $\un x=(x_1,\ldots,x_s)$ of elements in $\fm$ is said to be a dd-sequence of $M$ if for all $i\in \{1,\ldots,s\}$ and all $s$-tuples of positive integers $(n_1,\ldots,n_s)$, the sequence $x_1^{n_1},\ldots,x_i^{n_i}$ is a d-sequence of the module $M/(x_{i+1}^{n_{i+1}},\ldots,x_s^{n_s})M$.
Then the main result of this section is the following characterization of a \sop to be a dd-sequence.
\begin{theorem}\label{chara}
Let $M$ be a finitely generated $R$-module and $\un x =(x_1,\ldots,x_d)$ a system of parameters of $M$. Then $\un x$ is a dd-sequence if and only if 
$$\chi_1(\un x(\un n);M)=\sum_{i=0}^{p_1(M)}n_1\ldots n_ie(x_1,\ldots,x_i;(0:x_{i+1})_{M/(x_{i+2},\ldots,x_d)M})$$
for all $n_1,\ldots , n_d>0$.
\end{theorem}
\noindent
It follows from this theorem that every p-standard \sop  of $M$ is a dd-sequence. However, the converse statement is wrong. Observe that there always exist p-standard systems of parameters of $M$ provided the ring $R$ admits a dualizing complex (see \cite {c2}). Hence, in this case,  there exist also systems of parameters of $M$ which are dd-sequences.

The Section 4 is devoted to show the following theorem, which is one of the main results of the paper.
\begin{theorem}\label{main}
Let $M$ be a finitely generated $R$-module, $\un x=(x_1,\ldots,x_d)$ a system of parameters of $M$. Suppose that $\un x$ is a dd-sequence of $M$. Then the $k$-th Euler-Poincar\'{e} characteristic $\chi_k(\un x(\un n);M)$ is a polynomial in $\un n$ for all $k> 0$. Moreover, in this case we have
$$\chi_k(\un x(\un n);M)=\sum_{i=0}^{p_k(M)}n_1\ldots n_ie(x_1,\ldots,x_i;(0:x_{i+1})_{H_{k-1}(x_{i+2},\ldots,x_d;M)}),$$
where  we set $e(x_1,\ldots ,x_i;(0:x_{i+1})_{H_{k-1}(x_{i+2},\ldots ,x_d;M)})=\ell ((0:x_1)_{H_{k-1}(x_2,\ldots , x_d;M)})$ for $i=0$.
\end{theorem}

Let $k$ be an integer with $k\geq p_1(M)$ and $\un x=(x_1,\ldots,x_d)$  a system of parameters which is a dd-sequence. Then it follows from Theorem \ref {chara} that $p_1(M/(x_1,\ldots , x_k)M) \leq 0$. Therefore $M/(x_1,\ldots , x_k)M$ is a generalized Cohen-Macaulay module of dimension $d-k$. Thus $\ell (H_{\fm}^i(M/(x_1,\ldots , x_k)M))<\infty$ for all $i=0,\ldots , d-k-1$. 

As an application of  Theorem \ref {main}, we prove
in  Section 5 a result concerning with the length of local cohomology modules as follows.
\begin{theorem}\label {main4}
Let $M$ be a finitely generated $R$-module and $\un x$  a system of parameters of $M$. Assume that $\un x$ is a dd-sequence. Then $\ell(H^i_{\fm}(M/(x_1^{n_1},\ldots,x_k^{n_k})M))$ is a polynomial in $\un n$ for all $k\geq p_1(M)$ and $ i<d-k$.
\end{theorem}
The last Section is to  focus on sequentially Cohen-Macaulay modules which were introduced first for graded algebras by Stanley [22]. We show that there always exists \sop which is a dd-sequence on a sequentially \cm module and their  partial Euler-Poincar\'{e} characteristic are easy to compute by  multiplicities of  the dimension filtration of the module. Therefore the following result is the main result of this section.
\begin{theorem}\label {5}  Let $\un x =(x_1,\ldots,x_d)$ be a system of parameters of    a sequentially \cm  module $M$. Suppose that   
$$M_0\subset M_1\subset \cdots \subset M_{t-1}\subset M_t=M$$ 
 is the dimension filtration of $M$ with $\dim M_i =d_i$. Then it holds

\noindent
 (a) The following statements are equivalent:

i) $M_i\cap (x_{d_i+1},\ldots,x_d)M=0$ for all $i=0,1,\ldots , t-1$.

ii) $ (x_{d_i+1},\ldots,x_d)M_i=0$ for all $i=0,1,\ldots , t-1$.

iii) For all $n_1, \ldots, n_d\geq 1$,
$$\chi_1(\un x(\un n);M)=\sum_{i=0}^{t-1}n_1\ldots n_{d_i}e(x_1,\ldots,x_{d_i};M_i).$$

iv) $\un x$ is a dd-sequence of $M$.

\noindent
(b) Suppose that $\un x$ is satisfied one of the equivalent conditions above,  then 
$$\chi_k(\un x(\un n);M)=\sum_{d_i\leq d-k}a_in_1\ldots n_{d_i}e(x_1,\ldots,x_{d_i};M_i),$$
where 
$$a_{t-i}=\sum_{\substack{j_1+\ldots+j_{i-1}\leq k-1\\k\leq j_1+\ldots+j_i}}(-1)^{j_1+\cdots+j_i-k}\binom{d_t-d_{t-1}}{j_1}\cdots \binom{d_{i+1}-d_i}{j_{t-i}},$$
and we stipulate that $\binom{a}{b}=0$ if $a<b$.

\end{theorem}

%%%%%%%%%%%%%
%%%%%%%%%%%%%%%%%%%%%

\section{Partial Euler-Poincar\'{e} characteristics}
\markright{Euler-Poincar\'{e} characteristics}

Throughout this paper we always denote by $(R,\fm)$  a commutative Noetherian local  ring  with the maximal ideal $\fm$ and by $M$ a finitely generated $R$-module of dimension $d$.  Let $\un x =(x_1,\ldots,x_d)$ be a \sop of $M$. Consider the Koszul complex $K(\un x; M)$ of $M$ with respect to the \sop $\un x$. Let  $H_i(\un x;M)$ be   the $i$-th Koszul homology module of $K(\un x;M)$. Then   the {\it $k$-th  Euler-Poincar\'{e} characteristic} of $M$ with respect to   $\un x$ (or  the partial Euler-Poincar\'{e} characteristic of $K(\un x;M)$) is defined by 
$$\chi_k(\un x;M)=\sum_{i=k}^\infty(-1)^{i-k}\ell(H_i(\un x;M)).$$
It is well-known that $\chi_0(\un x; M)=e(\un x; M)$, $\chi_1(\un x; M)=\ell (M/\un xM)- e(\un x; M)$, $\chi_k(\un x;M)\geq 0$ for all $k$ and $\chi_k(\un x;M)=0$  for all $k>d$. The following lemma is a generalization of a result by M. Auslander and D. A. Buchsbaum (see \cite[Corollary 4.3]{ab}) and is often used  in the sequel.
\begin{lemma}{\rm \cite [ Corollary 2.2]{ck}}\label{ckh}
Let $\un x=(x_1,\ldots, x_d)$ be a \sop of $M$. Then 
$$\chi_k(\un x;M)=\sum_{i=0}^{d-k}e(x_1,\ldots,x_i;(0:x_{i+1})_{H_{k-1}(x_{i+2},\ldots,x_d;M)}).$$
\end{lemma}

Let now $\un n=(n_1,\ldots ,n_d)$ be  a d-tuple of positive integers. Put $\un x(\un n)=(x_1^{n_1},\ldots,x_d^{n_d})$ and consider $\chi_k(\un x(\un n);M)$ as a function in $\un n$.
 
In order to prove Theorem  \ref {main3} we need the following  auxiliary lemmata.
\begin{lemma} \label{l1} Let $\un x=(x_1,\ldots, x_d)$ be a \sop of $M$. Then the following statements are true for all $k\geq 0$.

(i) $\chi_k(\un x(\un n);M)\leq n_1\ldots n_d\chi_k(\un x;M)$ for all $\un n$.

(ii) $\chi_k(\un x(\un n);M)$ is an increasing function, i. e., for all $\un m =(m_1,\ldots , m_d)$, $n_1\leq m_1,\ldots , n_d\leq m_d$,
$$\chi_k(\un x(\un n);M)\leq \chi_k(\un x(\un m);M).$$

(iii) If $\chi_k(\un x(\un n);M)$ is a polynomial in $\un n$ then it is linear in each variable $n_i$,
 i. e.,  there exist integers  $\lambda_{i_1\ldots i_t}$ such that
$$\chi_k(\un x(\un n);M)=\sum_{t=0}^d \sum_{0<i_1<\ldots <i_t\leq d}\lambda_{i_1\ldots i_t}n_{i_1}\ldots n_{i_t},$$
where we set $\lambda_{i_1\ldots i_t}n_{i_1}\ldots n_{i_t}=\lambda $ the constant coefficient of the polynomial when $t=0$.
\end{lemma}
\begin {proof} The statements (i), (ii) were proved in \cite [Lemma 3.2] {ck} and (iii) follows easily from (i).
\end {proof}
\begin{remark}\label {r1}
 Garcia Roig and Kirby proved in \cite {rk} that the function $\chi_k(\un x(\un n);M)$  in general is not a polynomial. But, it was shown in \cite {roig} that for the case $n_1=n_2=\ldots =n_d =n$ the function $\chi_k(\un x( n);M)$, considered as a function in one variable $n$, is always bounded above by a polynomial in $n$ of degree at most $d-k$. Later, a more general result has been proved in \cite {ck}, which says  that the least degree of all polynomials in $\un n$ bounding above the function $\chi_k(\un x(\un n);M)$ is independent of the choice of the \sop $\un x$. Denote this new invariant by $p_k(M)$. Then we  have
$$p_0(M)=d> p_1(M)\geq \ldots \geq p_{d-1}(M)\geq p_d(M).$$
Especially, the invariant $p_1(M)$ was called in \cite {c2} the polynomial type of $M$ and denoted by $p(M)$.
\end {remark}
The next lemma easily follows from properties of linear polynomials. 
\begin {lemma}\label {l2}
Let $\varphi(n_1,\dots,n_d): \mathbb N^d \longrightarrow \mathbb N \cup \{ 0\} $ be a function defining on $\mathbb N^d$ ($\mathbb N$ is the set of all positive integers). Assume that  $\varphi(n_1,\dots,n_d)$ is a linear polynomial in each variable $n_i$ for all $i= 1,\ldots,d$.  Then $\varphi(n_1,\dots,n_d)$ is a polynomial in $\un n$.
\end{lemma}

\noindent
{\it Proof of Theorem \ref {main3}}.

\noindent
{\it Necessary condition}: Assume that $\chi_k(\un x(\un n);M)$ is a polynomial for all positive integers $n_1,\ldots, n_d\geq n_0 $. Since the function $\chi_k(\un x(\un n);M)$ does not depend on the order of the sequence $\un x$, we need only to check the necessary condition for the case $i=d$. Using Lemma \ref{ckh} to the sequence $x_d,x_{d-1},\ldots , x_1$  we have
$$\chi_k(\un x(\un n);M)=n_d\sum_{i=k}^{d-1}n_{i+1}\ldots n_{d-1}e(x_{i+1},\ldots ,x_d;(0:x_i^{n_i})_{H_{k-1}(x_1^{n_1},\ldots ,x_{i-1}^{n_{n-1}};M)})$$
$$ + \ell((0:x_d^{n_d})_{H_{k-1}(x_1^{n_1},\ldots ,x_{d-1}^{n_{d-1}};M)}).$$
Since $\chi_k(\un x(\un n);M)$ is a polynomial for all $n_i\geq n_0$ by the hypothesis, $\chi_k(\un x(\un n);M)$ is linear in $n_d$ by Lemma \ref{l1}, (iii). Therefore   
$\ell((0:x_d^{n_d})_{H_{k-1}(x_1^{n_1},\ldots ,x_{d-1}^{n_{d-1}};M)})$ 
is a polynomial in $n_d$, for all $n_i\geq n_0$.  
Thus 
$$(0:x_d^{n_d})_{H_{k-1}(x_1^{n_1},\ldots ,x_{d-1}^{n_{d-1}};M)} =(0:x_d^{n_0})_{H_{k-1}(x_1^{n_1},\ldots ,x_{d-1}^{n_{d-1}};M)},$$
for all $n_1,\ldots, n_d\geq n_0$.

\noindent
{\it Sufficient condition}: Without loss of generality we may assume that $n_0=1$. 
For each $i\in \{1,\ldots,d\}$, applying Lemma \ref{ckh} to the sequence  $x_i^{n_i},x_1^{n_1},\ldots,\widehat{x_i^{n_i}},$ $\ldots,x_d^{n_d}$, where $\widehat x$ indicates that $x$ is omitted there,  we obtain by the hypothesis that
\[ \begin{aligned}
\chi_k(\un x(\un n);M)&=n_i\phi_i(n_1,\ldots,\widehat{n_i},\ldots,n_d) + \ell((0:x_i^{n_i})_{H_{k-1}(x_1^{n_1},\ldots,\widehat{x_i^{n_i}},\ldots,x_d^{n_d};M)})\\
&=n_i\phi_i(n_1,\ldots,\widehat{n_i},\ldots,n_d) + \ell((0:x_i)_{H_{k-1}(x_1^{n_1},\ldots,\widehat{x_i^{n_i}},\ldots,x_d^{n_d};M)}),
\end{aligned} \]
where the function $\phi_i(n_1,\ldots,\widehat{n_i},\ldots,n_d)$ is  independent of $n_i$. Therefore, for all given (d-1)-tuple of positive integers $(n_1,\ldots,\widehat{n_i},\ldots,n_d) \in \mathbb N^{d-1}$, the function $\chi_k(\un x(\un n);M)$ is a linear polynomial  in $n_i$ for all $i=1,\ldots , d$. Hence, $\chi_k(\un x(\un n);M)$ is a polynomial in $\un n$ by Lemma \ref {l2}.
\qed
%%%%%%%%%%%%%%%
%%%%%%%%%%%%%%%%%%%%%%

\section{dd-Sequences}
\markright{dd-sequences}
First, we  recall the notion of d-sequence. 
\begin {definition}
 Let $\un x=(x_1,\ldots,x_s)$ be  a sequence of elements in the maximal ideal $\fm$. Then $\un x$ is called  a {\it d-sequence} of the $R$-module $M$ if  
$$(x_1,\ldots,x_{i-1})M:x_j=(x_1,\ldots,x_{i-1})M:x_ix_j$$
for all $i=1,\ldots,s$ and all $j\geq i$. 
Moreover, $\un x$ is called a strong d-sequence of $M$ if $\un x(\un n)=(x_1^{n_1},\ldots,x_s^{n_s})$ is a d-sequence for all $\un n=(n_1,\ldots , n_s) \in \mathbb{N}^s$.
\end {definition}

It should be mentioned that the notion of d-sequence was introduced  by C. Huneke \cite {ch} and it has become a useful tool in different topics of commutative algebra. Now we define a new notion of a sequence.
\begin{definition}
A sequence  $\un x=(x_1,\ldots,x_s)$  of elements in the maximal ideal $\fm$  is  called   a {\it dd-sequence} of $M$ if $\un x$ is a strong d-sequence and  the following inductive conditions are satisfied:  

(i) either $s=1$, or

(ii) $s>1$ and  the sequence $x_1,\ldots,x_{s-1}$ is a dd-sequence of $M/x_s^nM$ for all positive integers $n$.
\end{definition}

\begin{remark}\label{hah} (i) It is easy to see that an $M$-sequence is always a dd-sequence. An unconditioned strong d-sequence is also a dd-sequence. Therefore a part of \sop of a Buchsbaum module is a dd-sequence.

\noindent
(ii) If $x_1,\ldots,x_s$ is a dd-sequence of $M$ then $x_1^{n_1},\ldots,x_s^{n_s}$ is a dd-sequence of $M$ for all  $s$-tuples of positive integers $(n_1,\ldots,n_s)\in \mathbb N^s$.

\noindent
(iii) It is not difficult to show by the definition that a sequence $\un x=(x_1,\ldots,x_s)$ is a dd-sequence of $M$ if and only if  $x_1^{n_1},\ldots,x_i^{n_i}$ is a d-sequence of the module $M/(x_{i+1}^{n_{i+1}},\ldots,x_s^{n_s})M$ for all $i\in \{1,\ldots,s\}$ and all $s$-tuples of positive integers $(n_1,\ldots,n_s)\in \mathbb N^s$, where we set $M/(x_{i+1}^{n_{i+1}},\ldots,x_s^{n_s})M =M$ when $i=s$.
\end{remark}

Let $x_1,\ldots,x_s$ be a sequence of elements of $\fm$. For  convenience,  from now on we denote  the sequence $x_1,\ldots,x_{i-1},x_{i+1},\ldots,x_s$ by $x_1,\ldots,\widehat{x_i},\ldots,x_s$.
\begin{proposition}\label{link}
Let $x_1,\ldots,x_s$ be a dd-sequence of $M$. Then  $x_1,\ldots,\widehat{x_i},\ldots,x_s$ is a dd-sequence of $M/x_iM$  for  $i\in \{1,\ldots,s\}$. 
\end{proposition}
\begin{proof}
The proposition is proved by induction on $s$. The assertion follows immediately from the definition for cases $s=1,2$. Now assume that $s>2$. By virtue of Remark \ref{hah} we  have only to prove the assertion for $1<i<s$. Since $x_1,\ldots,x_{s-1}$ is  a dd-sequence of $M/x_sM$ by the definition, it follows from the inductive hypothesis that $x_1,\ldots,\widehat{x_i},\ldots,x_{s-1}$ is a dd-sequence of $ M/(x_i,x_s)M $. Therefore, in order to prove that $x_1,\ldots,\widehat{x_i},\ldots,x_s$ is  a dd-sequence of $M/x_iM$, we have  to verify that $x_1,\ldots,\widehat{x_i},\ldots,x_s$ is a strong d-sequence of $M/x_iM$. Moreover,  since $x_1^{n_1},\ldots, x_{s}^{n_{s}}$ is also a dd-sequence of $M$ for all  $(n_1,\ldots,n_{s})\in \mathbb N^{s}$ by Remark \ref {hah}, (ii), we need only to show that $x_1,\ldots,\widehat{x_i},\ldots,x_s$ is a  d-sequence of $M/x_iM$.
In fact, since $x_1,\ldots,x_{s-1}$ is a dd-sequence of $M/x_s^n$ for all $n>0$, we can prove by  using of Remark \ref {hah}, (iii) and Krull's Intersection Theorem that
 $x_1,\ldots,x_{s-1}$ is a dd-sequence of $M$. Therefore $x_1,\ldots,\widehat{x_i},\ldots,x_{s-1}$ is a d-sequence of $M/x_iM$ by inductive hypothesis. On the other hand, since 
$x_2,\ldots,x_s$ is a dd-sequence of $M/x_1M$,  $x_2,\ldots,\widehat{x_i},\ldots,x_s$ is a dd-sequence of $M/(x_1,x_i)M$ by the inductive hypothesis. 
Combine these facts  and keep in mind of the definition of a d-sequence,  it remains to check that $0:_{M/x_iM}x_s=0:_{M/x_iM}x_1x_s$ or equivalent, $x_iM:x_s= x_iM:x_1x_s$. First, we show that  
$x_iM:x_1\subseteq x_iM:x_s$. To do this 
let $a$ be an arbitrary element of $x_iM:x_1$. Then, there exists $b\in M$ such that $x_1a=x_ib$, hence $b\in x_1M:x_i \subseteq x_1M:x_s$. Thus there exists an element $a^\prime \in M$ such that $x_sb=x_1a^\prime$. Therefore, $x_1x_sa=x_ix_sb=x_1x_ia^\prime$. It follows that $x_sa-x_ia^\prime \in 0:_Mx_1\subseteq 0:_Mx_s$,  and hence $x_s(x_sa-x_ia^\prime )=0$. This leads to $x_s^2a=x_sx_ia^\prime \in x_iM$. So we get $a\in x_iM:x_s^2$.  Since $\un x$ is a dd-sequence of $M$,  $x_iM:x_s^2=x_iM:x_s$. Therefore  $x_iM:x_1\subseteq x_iM:x_s$. 
Finally, we obtain
\[ \begin{aligned}
x_iM:x_1x_s&=(x_iM:x_1):x_s\\
&\subseteq (x_iM:x_s):x_s
=x_iM:x_s^2=x_iM:x_s
\end{aligned} \]
and the proof  is  complete.
\end{proof}

\begin{corollary}\label{hq}
Let  $x_1,\ldots,x_s$ be  a dd-sequence of $M$. Then for all $i\in \{ 1, \ldots , s\}$, the sequence $x_1,\ldots,\widehat{x_i},\ldots,x_s$ is also a dd-sequence of $M$.
\end{corollary}
\begin{proof}
Since $x_1,\ldots,x_s$ is a dd-sequence,  $ x_1,\ldots,x_i^n,\ldots,x_s$ is a dd-sequence of $M$ for all $n>0$. Hence, by Proposition \ref{link}, $x_1,\ldots,\widehat{x_i},\ldots,x_s$ is also a dd-sequence of $M/x_i^nM$. Therefore the corollary follows easily from Remark \ref {hah}, (iii) and  Krull's  Intersection Theorem. 
\end{proof}
Now we  prove the main result of this Section.

\noindent
{\it Proof of Theorem \ref {chara}}.

\noindent
{\it Necessary condition}: Recall that if $\chi_1(\un x(\un n);M)$ is a polynomial in $\un n$, then this polynomial must have the degree $p_1(M)\leq d-1$. Therefore the statement is proved if we can show that
$$\chi_1(\un x(\un n);M)=\sum_{i=0}^{d-1}n_1\ldots n_ie(x_1,\ldots,x_i;(0:x_{i+1})_{M/(x_{i+2},\ldots,x_d)M})$$
for all $n_1,\ldots , n_d>0$. We proceed by induction on $d$. 
 If $d=1$ there is nothing to do. Assume that $d>1$. By Lemma  \ref{ckh} we have 
$$\chi_1(\un x(\un n);M)=\sum_{i=1}^dn_{i+1}\ldots n_de(x_{i+1},\ldots,x_d;(0:x_i^{n_i})_{M/(x_1^{n_1},\ldots,x_{i-1}^{n_{i-1}})M}).$$
Since $\un x$ is a strong d-sequence of $M$,
\[ \begin{aligned}
(0:x_i^{n_i})_{M/(x_1^{n_1},\ldots,x_{i-1}^{n_{i-1}})M}&\subseteq (0:x_d^{n_d})_{M/(x_1^{n_1},\ldots,x_{i-1}^{n_{i-1}})M}\\
&=(0:x_d)_{M/(x_1^{n_1},\ldots,x_{i-1}^{n_{i-1}})M}
\end{aligned} \]
for all $i\leq d$.
Thus 
$$e(x_{i+1},\ldots,x_d;(0:x_i^{n_i})_{M/(x_1^{n_1},\ldots,x_{i-1}^{n_{i-1}})M})=0, \ i<d,$$
and 
$$\ell((0:x_d^{n_d})_{M/(x_1^{n_1},\ldots,x_{d-1}^{n_{d-1}})M})=\ell((0:x_d)_{M/(x_1^{n_1},\ldots,x_{d-1}^{n_{d-1}})M}).$$
Therefore the function $\chi_1(\un x(\un n);M)$ is independent of $n_d$. Put $M^\prime =M/x_dM$ and $\un x'(\un n')= (x_1^{n_1},\ldots,x_{d-1}^{n_{d-1}})$.  Then we have
\[ \begin{aligned} 
\chi_1(\un x(\un n);M)&=\chi_1(\un x'(\un n'),x_d;M)\\
&=\chi_1(\un x^\prime(\un n^\prime);M^\prime)+n_1\ldots n_{d-1}e(x_1,\ldots,x_{d-1};0:_Mx_d).
\end{aligned} \]
It follows by the inductive hypothesis that
\[ \begin{aligned} 
\chi_1(\un x'(\un n');M^\prime)&=\sum_{i=0}^{d-2}n_1\ldots n_ie(x_1,\ldots,x_i;(0:x_{i+1})_{M^\prime/(x_{i+2},\ldots,x_{d-1})M^\prime})\\
&=\sum_{i=0}^{d-2}n_1\ldots n_ie(x_1,\ldots,x_i;(0:x_{i+1})_{M/(x_{i+2},\ldots,x_{d-1},x_d)M}).
\end{aligned} \]
Therefore 
$$\chi_1(\un x(\un n);M)=\sum_{i=0}^{d-1}n_1\ldots n_ie(x_1,\ldots,x_i;(0:x_{i+1})_{M/(x_{i+2},\ldots,x_d)M})$$
as required.

\noindent
{\it Sufficient condition}: We also prove the statement by induction on $d$. The case $d=1$ is trivial. Suppose that $d>1$.  Set $\overline M=M/x_d^{n_d}M $  for $n_d\in \mathbb N$. Then  we have
\[ \begin{aligned}
\chi_1(\un x^\prime(\un n^\prime);\overline M)&=\chi_1(\un x(\un n);M)-n_1\ldots n_{d-1}e(x_1,\ldots,x_{d-1};0:_Mx_d^{n_d})\\
&=\sum_{i=0}^{d-2}n_1\ldots n_ie(x_1,\ldots,x_i;(0:x_{i+1})_{M/(x_{i+2},\ldots,x_{d-1},x_d)M})\\
&+ n_1\ldots n_{d-1}(e(x_1,\ldots,x_{d-1};0:_Mx_d)- e(x_1,\ldots,x_{d-1};0:_Mx_d^{n_d})).
\end{aligned} \]
Since $\chi_1(\un x^\prime(\un n^\prime);\overline M)$ forms a polynomial in $\un n$, $$\text {deg }(\chi_1(\un x^\prime(\un n^\prime);\overline M))=p_1(\overline M)\leq \dim \overline M-1=d-2.$$ Therefore
$$e(x_1,\ldots,x_{d-1};0:_Mx_d)- e(x_1,\ldots,x_{d-1};0:_Mx_d^{n_d})=0.$$
Hence 
$$\chi_1(\un x^\prime(\un n^\prime);\overline M)=\sum_{i=0}^{d-2}n_1\ldots n_ie(x_1,\ldots,x_i;(0:x_{i+1})_{\overline M/(x_{i+2},\ldots,x_{d-1})\overline M}).$$
It follows by inductive hypothesis that $x_1,\ldots,x_{d-1}$ is a dd-sequence of $\overline M$. Therefore we have only to prove that $\un x$ is a strong d-sequence of $M$.  Since $x_1,\ldots,x_{d-1}$ is a strong d-sequence of $\overline M$, 
$$(x_d^{n_d},x_1^{n_1},\ldots,x_i^{n_i})M:x_j^{n_j}=(x_d^{n_d},x_1^{n_1},\ldots,x_i^{n_i})M:x_{i+1}^{n_{i+1}}x_j^{n_j}$$
for all $0\leq i <j \leq d-1$.
Hence, by Krull's Intersection Theorem,
\[ \begin{aligned}
(x_1^{n_1},\ldots,x_i^{n_i})M:x_j^{n_j}
&=\bigcap_{n_d}(x_d^{n_d},x_1^{n_1},\ldots,x_i^{n_i})M:x_j^{n_j}\\
&=\bigcap_{n_d}(x_d^{n_d},x_1^{n_1},\ldots,x_i^{n_i})M:x_{i+1}^{n_{i+1}}x_j^{n_j}\\
&=(x_1^{n_1},\ldots,x_i^{n_i})M:x_{i+1}^{n_{i+1}}x_j^{n_j}.
\end{aligned} \]
So in order to show that $\un x$ is a strong d-sequence of $M$. It remains to check that $$(x_1^{n_1},\ldots,x_i^{n_i})M:x_d^{n_d}=(x_1^{n_1},\ldots,x_i^{n_i})M:x_{i+1}^{n_{i+1}}x_d^{n_d},\hbox{ for all } 0\leq i<d.$$
In fact, applying Lemma \ref{ckh} to the sequence $x_d^{n_d},x_{i+1}^{n_{i+1}},x_1^{n_1},\ldots,\widehat{x_{i+1}^{n_{i+1}}},\ldots,x_{d-1}^{n_{d-1}}$, we obtain for all $i=0,\ldots , d-1$ that
\[ \begin{aligned}
\chi_1(\un x(\un n);M)
&=\sum_{\substack{j=0\\j\not= i}}^{d-2}
e(x_d^{n_d},x_{i+1}^{n_{i+1}},x_1^{n_1},\ldots,x_j^{n_j};(0:x_{j+1}^{n_{j+1}})_{M/(x_{j+2}^{n_{j+2}},\ldots,\widehat{x_{i+1}^{n_{i+1}}},\ldots,x_{d-1}^{n_{d-1}})M})\\
&\hspace{0.3cm}+n_de(x_d;(0:x_{i+1}^{n_{i+1}})_{M/(x_1^{n_1},\ldots,\widehat{x_{i+1}^{n_{i+1}}},\ldots,x_{d-1}^{n_{d-1}})M})\\
&\hspace{0.3cm}+\ell((0:x_d^{n_d})_{M/(x_{i+1}^{n_{i+1}},x_1^{n_1},\ldots,x_{d-1}^{n_{d-1}})M}).
\end{aligned} \]
Moreover, since $\chi_1(\un x(\un n);M)$ does not depend on $n_d$, we deduce that $$e(x_d;(0:x_{i+1}^{n_{i+1}})_{M/(x_1^{n_1},\ldots,\widehat{x_{i+1}^{n_{i+1}}},\ldots,x_{d-1}^{n_{d-1}})M})=0.$$
Thus there exists a positive integer  $n$ such that $$x_d^n(0:x_{i+1}^{n_{i+1}})_{M/(x_1^{n_1},\ldots,\widehat{x_{i+1}^{n_{i+1}}},\ldots,x_{d-1}^{n_{d-1}})M}=0.$$ 
So
$$(0:x_{i+1}^{n_{i+1}})_{M/(x_1^{n_1},\ldots,\widehat{x_{i+1}^{n_{i+1}}},\ldots,x_{d-1}^{n_{d-1}})M}\subseteq (0:x_d^n)_{M/(x_1^{n_1},\ldots,\widehat{x_{i+1}^{n_{i+1}}},\ldots,x_{d-1}^{n_{d-1}})M}.$$
On the other hand, since $\chi_1(\un x(\un n);M)$ is a polynomial, it follows  by \cite[Theorem 1]{c1} that
$$(0:x_d^{n+n_d})_{M/(x_1^{n_1},\ldots,\widehat{x_{i+1}^{n_{i+1}}},\ldots,x_{d-1}^{n_{d-1}})M}\\
=(0:x_d^{n_d})_{M/(x_1^{n_1},\ldots,\widehat{x_{i+1}^{n_{i+1}}},\ldots,x_{d-1}^{n_{d-1}})M}$$
for all $n_d$. Therefore
\[ \begin{aligned}
(0:x_{i+1}^{n_{i+1}}x_d^{n_d})_{M/(x_1^{n_1},\ldots,\widehat{x_{i+1}^{n_{i+1}}},\ldots,x_{d-1}^{n_{d-1}})M}&\subseteq (0:x_d^{n+n_d})_{M/(x_1^{n_1},\ldots,\widehat{x_{i+1}^{n_{i+1}}},\ldots,x_{d-1}^{n_{d-1}})M}\\
&=(0:x_d^{n_d})_{M/(x_1^{n_1},\ldots,\widehat{x_{i+1}^{n_{i+1}}},\ldots,x_{d-1}^{n_{d-1}})M}.
\end{aligned} \]
Thus $$(x_1^{n_1},\ldots,\widehat{x_{i+1}^{n_{i+1}}},\ldots,x_{d-1}^{n_{d-1}})M:x_d^{n_d}=(x_1^{n_1},\ldots,\widehat{x_{i+1}^{n_{i+1}}},\ldots,x_{d-1}^{n_{d-1}})M:x_{i+1}^{n_{i+1}}x_d^{n_d}.$$
Now,  applying Krull's  Intersection Theorem again we get
\[ \begin{aligned}
(x_1^{n_1},\ldots,x_i^{n_i})M:x_d^{n_d}&=\bigcap_{n_{i+2},\ldots,n_{d-1}}(x_1^{n_1},\ldots,\widehat{x_{i+1}^{n_{i+1}}},\ldots,x_{d-1}^{n_{d-1}})M:x_d^{n_d}\\
&=\bigcap_{n_{i+2},\ldots,n_{d-1}}(x_1^{n_1},\ldots,\widehat{x_{i+1}^{n_{i+1}}},\ldots,x_{d-1}^{n_{d-1}})M:x_{i+1}^{n_{i+1}}x_d^{n_d}\\
&=(x_1^{n_1},\ldots,x_i^{n_i})M:x_{i+1}^{n_{i+1}}x_d^{n_d}
\end{aligned} \]
for all $i=0,\ldots ,d-1$. So $\un x$ is a strong d-sequence and the proof of Theorem 1.2 is complete.\qed

\medskip
As an application of Theorem \ref {chara}  we get the following useful criterion to check whether a \sop is a dd-sequence.
\begin{corollary}
A system of parameters $\un x$ of $M$  is a dd-sequence of $M$ if and only if there exist integers $a_0, a_1,\ldots , a_{d-1}$ such that
$$\chi_1(\un x(\un n);M)=\sum_{i=0}^{d-1} n_1\ldots n_ia_i.$$ 
Moreover, in this case it holds
$$a_i=e(x_1,\ldots, x_i;(0:x_{i+1})_{M/(x_{i+2},\dots,x_d)M})$$
for $i=0,1,\ldots , p_1(M)$ and $a_{p_1(M)+1}=\ldots =a_{d-1}=0$.
\end{corollary}
\begin{proof}
The necessary condition follows immediately Theorem \ref{chara}. We proceed by induction on $d$ that
$$a_i=e(x_1,\ldots, x_i;(0:x_{i+1})_{M/(x_{i+2},\dots,x_d)M})$$
for all $i=0,1, \ldots , d-1,$ and the corollary follows by Theorem \ref {chara}. In fact,
 the case $d=1$ is obvious.  Suppose that $d>1$. 
We set $M_d=M/x_dM$ and $\un x'(\un n')= (x_1^{n_1},\ldots,x_{d-1}^{n_{d-1}})$. Then 
\[ \begin{aligned}
\chi_1(\un x'(\un n');M_d)
&=\chi_1(\un x(\un n);M)-e(x_1^{n_1},\ldots, x_{d-1}^{n_{d-1}};0:_Mx_d)\\
&=\sum_{i=0}^{d-2}n_1\ldots n_ia_i+n_1\ldots n_{d-1}(a_{d-1}-e(x_1,\ldots, x_{d-1};0:_Mx_d)).
\end{aligned} \]
Since the function $\chi_1(x_1^{n_1},\ldots,x_{d-1}^{n_{d-1}};M_d)$ is increasing and bounded above by a polynomial in $n_1,\ldots , n_{d-1}$ of degree $d-2$, it implies that 
$$a_{d-1}=e(x_1,\ldots, x_{d-1};0:_Mx_d).$$  
Therefore 
$$\chi_1(x_1^{n_1},\ldots,x_{d-1}^{n_{d-1}};M_d)
=\sum_{i=0}^{d-2}n_1\ldots n_ia_i.$$
So we get   by the inductive hypothesis that 
$$a_i=e(x_1,\ldots, x_i;(0:x_{i+1})_{M/(x_{i+2},\dots,x_d)M}),$$
for $i=0,\ldots, d-2$. 
\end{proof}
The following corollary is an immediate consequence of Theorem \ref {chara}.
\begin{corollary}
Let $\un x=(x_1,\ldots , x_d)$ be a \sop of $M$ and $k=p_1(M)$.   Suppose that $\un x$ is a dd-sequence of $M$. Then, for any  permutation $\pi$ of the set $\{ k+1,\ldots,d\}$, $x_1,\ldots,x_k,x_{\pi(k+1)},\ldots,x_{\pi(d)}$  is also a dd-sequence .
\end{corollary}
\begin {remark}
(i) Recall that an $R$-module $M$ is said to be a generalized Cohen-Macaulay module, if $ \underset { \un x } { \sup} \{ \chi_1(\un x;M)\}=I(M)<\infty$, where $\un x$ runs through all systems of parameters of $M$. It was shown in \cite {tr2} that $M$ is a generalized Cohen-Macaulay module if and only if there exists a \sop  $\un x$ of $M$ such that $ \chi_1(\un x;M)=I(M)$. That \sop  was called a standard \sop of $M$ and plays an important role in the theory of generalized Cohen-Macaulay modules. Suppose now that $M$ is a generalized Cohen-Macaulay module and $\un x $  a \sop of $M$. Then, it follows from Theorem \ref {chara} that $\un x$ is a standard \sop if and only if $\un x $ is a dd-sequence. Moreover, an $R$-module $M$ is called a Buchsbaum module if every \sop of $M$ is a standard system of parameters. Therefore we get the following characterization of Buchsbaum modules by dd-sequences as follows: $M$ is a Buchsbaum module  if and only if every \sop of $M$ is a dd-sequence.

\noindent
(ii) Let $\fa_i(M)=\text {Ann}_R(H_\fm^i(M))$ be  the annihilator ideal of the $i$-th local cohomology module  $H_\fm^i(M)$  of $M$ and $\fa(M)=\fa_0(M)\ldots\fa_{d-1}(M),$ where $d=\dim M$. Then, a \sop $\un x$ of $M$ was called a {\it p-standard system of parameters} if $x_i\in \fa(M/(x_{i+1},\ldots,x_d)M)$ for $i=1,\ldots,d$. This kind of \sop was introduced and investigated in \cite {c1', c2, c3} and has been used for solving the problem of Macaulayfication by recent works of T. Kawasaki \cite {ta1, ta2}. It was proved in \cite {c3} that if $\un x$ is a p-standard \sop  of $M$, then 
$$\chi_1(\un x(\un n);M)=\sum_{i=0}^{p_1(M)}n_1\ldots n_ie(x_1,\ldots,x_i;(0:x_{i+1})_{M/(x_{i+2},\ldots,x_d)M}).$$
Therefore any p-standard \sop is a dd-sequence.
Unfortunately, the converse is wrong. It means that a system of parameters, which is a dd-sequence,  in general is  not a p-standard \sop (see Example \ref {exp2}). However, the following consequence gives us close relations between them. 
\end {remark}

\begin{corollary}\label{a(M)}
Let $\un x=(x_1,\ldots , x_d)$ be a \sop of  $M$. If $\un x$ is a p-standard  \sop then it is also a dd-sequence. Conversely, suppose that $\un x$  is  a dd-sequence. Then $x_1^{n_1},\ldots,x_d^{n_d}$ is a p-standard system of parameters for all $n_i\geq i,\ i=1,\ldots,d$.
\end{corollary}
\begin{proof}The first statement is already shown in the remark above. Since a dd-sequence is a strong d-sequence, the second one is just an immediate consequence of Lemma 2.9 in \cite {c3}, which says that
$$x_jH_\fm^i(M/(x_1^{n_1},\ldots,x_h^{n_h})M)=0$$
for all $j=1,\ldots,d,\ h+i<j$ and $n_1,\ldots,n_h>0$, provided the system of parameters  $\un x=(x_1,\ldots,x_d)$ is a strong d-sequence. It should be noticed that this result was stated in \cite {c3} for p-standard systems of parameters, but in its proof one needs only the condition that the \sop  $\un x$ is a strong d-sequence.  
\end{proof}

 As we have seen,  a p-standard \sop is a dd-sequence and a dd-sequence is always a strong d-sequence.  However, the converse statements are wrong in general. Below we present some examples to clarify these three notions.
\begin{example}
Let $R= k[[X,Y]]$ be the ring  of all formal power series of two indeterminates $X, Y$ over a field $k$ with the maximal ideal $\fm$. Consider  $M=\fm^2$ as an $R$-module. Then $\dim(M)=2$ and $X,Y^2$ form  a \sop of $M$.  It is easy to check that $0:_MX^mY^{2n}=0=0:_MY^{2n}$ and $X^mM:Y^{2n}=X^mM:Y^2$ for all positive integers $m,n$. Therefore $X,Y^2$ is a strong d-sequence of $M$. On the other hand, we can show that
\begin{equation}\notag  
Y^2M:X^m=
\begin{cases}
(XY^2,Y^3), &\text{ if } m=1,\\
(Y^2),&\text{ if } m>1.
\end{cases}   
\end{equation}
Hence $X, Y^2$ is not a dd-sequence of $M$. Besides, we can also compute the $1$-st Euler-Poincar\'{e} characteristic and get
\begin{equation}\notag  
\chi_1(X^m,Y^{2n};M)=
\begin{cases}
2, &\text{ if } m=1,\\
3,&\text{ if } m>1.
\end{cases}   
\end{equation}
\end{example}
\begin{example}\label {exp2}
Let $R=k[[X_1,\ldots,X_{d+1}]], \ (d>1)$ be the ring of all formal power series of $(d+1)$-indeterminates $X_1, \ldots , X_{d+1}$ over a field $k$ with the maximal ideal $\fm =(X_1,\ldots , X_{d+1})$. Let  $M$ denote the $R$-module $R/I$ where $I=(X_{d+1}^{d+1},X_1X_{d+1}^d,X_2X_{d+1}^{d-1},\ldots,X_dX_{d+1})$. It is easy to see that $\dim(M)=d$  and $X_1,\ldots,X_d$ is a \sop  of $M$. Then by  a simple computation we get for all  $n_1,\ldots,n_d\in \mathbb N^d$, 
$$\ell(M/(X_1^{n_1},\ldots,X_d^{n_d})M)=\sum_{i=0}^dn_1\ldots n_i,$$
where we set $n_1\ldots n_i=1$ if $i=0$. Hence 
$$\chi_1(X_1^{n_1},\ldots,X_d^{n_d};M)=\sum_{i=0}^{d-1}n_1\ldots n_i.$$
Therefore the \sop  $X_1,\ldots,X_d$ is a dd-sequence of $M$ with $p_1(M)=d-1$. 
On the other hand, we can check that 
$$H_\fm^0(M)\cong (X_{d+1}^d,X_2X_{d+1}^{d-1},\ldots,X_dX_{d+1})/I.$$
 Hence 
$\fa_0(M)=\mathrm{Ann}_R(H_\fm^0(M))=\fm $.
 Moreover, since $R$ is  a regular local ring,  we have $\dim(R/\fa(M))=p_1(M)=d-1>0$ by \cite[Theorem 1.2]{c1'}. It follows that there exists at least an $i\in \{1,\ldots,d-1\}$ such that $\fa_i(M)\subseteq \fm$. Therefore $\fa(M)\subseteq \fa_0(M)\fa_i(M)\subseteq \fm^2$. 
Since $X_d\notin \fm^2$, $X_d\notin \fa(M)$. Thus $X_1,\ldots,X_d$ is not a p-standard $\sop$ of $M$.
\end{example}

%%%%%%%%%%%
%%%%%%%%%%%

\section{Partial Euler-Poincar\'{e} characteristics with respect to a dd-sequence}
\markright{Higher partial Euler-Poincar\'{e} characteristics}
Keep all notations in the previous sections. We  begin this Section with the following key lemma.
\begin{lemma}\label{lm}
Suppose that $\un x=(x_1,\ldots,x_s)$ is a dd-sequence of $M$. Then 
$$(0:x_j)_{H_k(x_1,\ldots,x_{i-1}, x_{j+1},\ldots,x_s;M)}=(0:x_ix_j)_{H_k(x_1,\ldots,x_{i-1}, x_{j+1},\ldots,x_s;M)}$$
for all $1\leq i\leq j\leq s$ and $0\leq k\leq s$. 
\end{lemma}
\begin{proof} Since the lemma is trivial for the case $s=1$, we assume that $s\geq 2$.  Moreover, since   $x_1,\ldots,x_{i-1},x_i,x_j, x_{j+1},\ldots,x_s$ is  a dd-sequence of $M$ by the hypothesis and Corollary \ref{hq},  it  suffices to prove the lemma for only  two cases $j=i+1$ and $j=i$. We will do it for the case $j=i+1$ and the proof for the other case  is  quite  similar. In fact, we need to show that
\begin {equation}\label {1}
(0:x_{i+1})_{H_k(x_1,\ldots,x_{i-1},x_{i+2},\ldots,x_s;M)}=(0:x_ix_{i+1})_{H_k(x_1,\ldots,x_{i-1},x_{i+2},\ldots,x_s;M)}. \tag {*}
\end {equation}
Set  $\un y=(y_1,\ldots y_s)$, where
\begin{equation}\notag  
y_j=
\begin{cases}
x_j &\text{ if } j<i,\\
x_{j+2}&\text{ if } i\leq j\leq s-2,\\
x_i&\text{ if } j=s-1,\\
x_{i+1}&\text{ if } j=s.
\end{cases}   
\end{equation}
Then the above equality (*)  can be rewritten as follows
\begin {equation} \label {2}
(0:y_s)_{H_k(y_1,\ldots,y_{s-2};M)}=(0:y_{s-1}y_s)_{H_k(y_1,\ldots,y_{s-2};M)}.\tag {**}
\end {equation}
Denote by $\varphi_k$ the $k$-th differential of the Koszul complex $K(y_1,\ldots,y_{s-2};M)$.
Since
$H_k(y_1,\ldots,y_{s-2};M)=\ker (\varphi_k)/\mathrm {im }(\varphi_{k+1})$, the equality (**) and hence the lemma are proved if we can verify the following claim.

\medskip
\noindent
 {\bf Claim.} {\it With all notations above it holds} 
$$ \mathrm{im}(\varphi_{k+1}):_{\ker (\varphi_k)}y_s=\mathrm{im}(\varphi_{k+1}):_{\ker(\varphi_k)}y_sy_{s-1}.$$ 
{\it for all $s\geq 2$ and} $k=0,1,\ldots , s$.
\medskip

\noindent
{\it Proof of the Claim}.
 Recall that any element $a\in K_k(y_1,\ldots,y_{s-2};M)\cong M^{\binom{s-2}{k}}$  can be expressed as $a=(a_{i_1,\ldots,i_k})_{1\leq i_1<i_2<\ldots<i_k\leq s-2 }\in M^{\binom{s-2}{k}}$. Extend  the  $ a_{i_1,\ldots,i_k}$'s  to  an  alternating function of the indices, i. e. such that $ a_{\ldots,i,\ldots,i,\ldots}=0$  and $ a_{\ldots,i,\ldots,j,\ldots}=-a_{\ldots,j,\ldots,i,\ldots}$ for all $i\not= j$. For convenience,  we say in this case that the elements  $ a_{i_1,\ldots,i_k}\in M$ are of alternating indices.
 Then we have
$$ \ker(\varphi_k)=\{(a_{i_1,\ldots,i_k})_{1\leq i_1<i_2<\ldots<i_k\leq s-2 }\in M^{\binom{s-2}{k}}| \sum_{j=1}^{s-2}y_ja_{j,i_1,\ldots,i_{k-1}}=0\},$$
\begin{multline*}
\mathrm {im} (\varphi_{k+1})=\{ (b_{i_1,\ldots,i_k})_{1\leq i_1<i_2<\ldots<i_k\leq s-2}\in M^{\binom{s-2}{k}}| \\ \exists  (a_{i_1,\ldots,i_{k+1}})_{1\leq i_1< \ldots <i_{k+1}\leq s-2}\in M^{\binom{s-2}{k+1}}:
 b_{i_1,\ldots,i_k}=\sum_{j= 1}^{s-2}y_ja_{j,i_1,\ldots,i_k}\}.
\end{multline*}
 We prove the claim by induction on $s$. The case $s= 2$ is straightforward. Now assume that $s >2$. 
Let $(a_{i_1,\ldots,i_k})_{1\leq i_1<\ldots <i_k\leq s-2}$ be  an arbitrary element
 of $ \mathrm{im}(\varphi_{k+1}):_{\ker(\varphi_k)}y_{s-1}y_s$. 
Then, it follows from  the above descriptions of $ \ker(\varphi_k)$ and $\mathrm {im} (\varphi_{k+1})$ that
$$(a_{i_1,\ldots,i_k})_{1\leq i_1<\ldots <i_k\leq s-2}\in  \mathrm{im}(\varphi_{k+1}):_{\ker(\varphi_k)}y_s$$
if   there exist elements $b_{j,i_1,\ldots,i_k} \in M$ of alternating indices such that 
\begin {equation}\label {3}
y_sa_{i_1,\ldots,i_k}=\sum_{j=1}^{s-2}y_jb_{j,i_1,\ldots,i_k}.\tag {***}
\end {equation}
Therefore  the claim is proved if we can show  the existence of $b_{j,i_1,\ldots,i_k}$ and the equality (***).
To do it, we consider the following cases.

\noindent
{\it Case 1}: $i_1>1$. 
Put $\overline M=M/y_1M$ and ${\overline a}_{i_1,\ldots,i_k}$  the image of $a_{i_1,\ldots,i_k}$ in $\overline M$. Denote by  $\varphi_k^\prime$  the $k$-th  differential of the Koszul complex $K(y_2,\ldots,y_{s-2};\overline M)$. Then it follows from
 Proposition \ref{link} and the inductive hypothesis that 
$$({\overline a}_{i_1,\ldots,i_k})_{ i_1<\ldots <i_k\leq s-2}
\in \mathrm{im}(\varphi_{k+1}^\prime):_{\ker (\varphi_k^\prime)}y_{s-1}y_s =\mathrm{im}(\varphi_{k+1}^\prime):_{\ker(\varphi_k^\prime)}y_s.$$
 So there exist elements $b_{j,i_1,\ldots,i_k} \in M$ of alternating indices  such that $$y_sa_{i_1,\ldots,i_k}=\sum_{j=1}^{s-2}y_jb_{j,i_1,\ldots,i_k}.$$

\noindent
{\it Case 2}: $i_1=1$. We have to prove in this case that 
$$y_sa_{1,i_2, \ldots,i_k}=\sum_{j=2}^{s-2}y_jb_{j,1,i_2, \ldots,i_k},$$
where $b_{j,1,\ldots,i_k}$'s are of alternating indices and just determined in the case 1.
Since $x_1,\ldots,x_s$ is a dd-sequence of $M$, we can apply  Proposition \ref {link} several times to imply that the sequence $x_i,x_{i+1}$ is a dd-sequence of the module $M/(x_1,\ldots,x_{i-1},x_{i+2},\ldots,x_s)M$. Therefore $y_{s-1},y_s$ is a d-sequence with respect to $M/(y_1,\ldots,y_{s-2})M$.  Hence, from the hypothesis
$$y_{s-1}y_sa_{1,i_2,\ldots,i_k}=\sum_{j= 2}^{s-2}y_ja_{j,1,i_2,\ldots,i_k},$$
it follows that
$$a_{1,i_2,\ldots,i_k} \in (\sum_{j=2}^{s-2}y_jR)M:y_{s-1}y_s=(\sum_{j=2}^{s-2}y_jR)M:y_s.$$
Let $t$ be a positive number such that  
$$y_sa_{1,i_2,\ldots,i_k}=\sum_{j=2}^{ t-1} y_jb_{j,1,i_2,\ldots,i_k}+\sum_{j=t}^{s-2}y_jc_j,$$
where $c_t,\ldots , c_{s-2}\in M$ . Then $t\geq 2$ and the first summand of the right term is $0$  if $t=2$.
 We put $I=(y_{t+1},\ldots y_{s-2})R$,  $M' =M/IM$ and denote by $a'$ the image of the element $a\in M$ in $M'$. We get 
$$y_sa'_{1,i_2,\ldots,i_k}=\sum_{j=2}^{t-1} y_jb'_{j,1,i_2,\ldots,i_k}+y_tc'_t$$
and by the case 1 that
$$y_sa'_{i_1,\ldots,i_
k}=\sum_{j= 1}^t y_jb'_{j,i_1,\ldots,i_k},\ (i_1>1).$$
Therefore (note that $b_{i,j,i_2,\ldots,i_k} \text{ are of alternating indices}$)
\[
\begin{aligned}
y_s\sum_{j=1}^ty_ja'_{j,i_2,\ldots,i_k}
&=\sum_{j=1}^ty_j(y_sa'_{j,i_2,\ldots,i_k})\\
&=\sum_{j=2}^ty_j(y_sa'_{j,i_2,\ldots,i_k})+y_1(y_sa'_{1,i_2,\ldots,i_k})\\
&=\sum_{j=2}^t
\sum_{i=1}^ty_jy_ib'_{i,j,i_2,\ldots,i_k}+\sum_{j=2}^{t-1} y_1y_jb'_{j,1,i_2,\ldots,i_k}+y_1y_tc'_t\\
&=y_1y_t(b'_{1,t,i_2,\ldots,i_k}+c'_t)\\
&=y_1y_t(c'_t-b'_{t,1,i_2,\ldots,i_k}).
\end{aligned}
\]
On the other hand, since $(a_{j,i_2,\ldots,i_k})_{j,i_2,\ldots,i_k}\in \ker (\varphi_k)$, 
$$\sum_{j=1}^{s-2}y_ja_{j,i_2,\ldots,i_k}=0,\ \hbox{ thus }\sum_{j=1}^ty_ja'_{j,i_2,\ldots,i_k}=0.$$
Hence $$0=y_s\sum_{j=1}^ty_ja'_{j,i_2,\ldots,i_k}=y_1y_t(c'_t-b'_{t,1,i_2,\ldots,i_k}).$$ 
Moreover, since $\un x$ is a dd-sequence of $M$, it is a strong d-sequence. We can show by applying  Krull's Intersection Theorem that $(x_1,\ldots ,x_{i-1},x_{i+1},\ldots , x_{s-2})=(y_1,\ldots , y_{s-2})$ is also a strong d-sequence.  Hence $y_1,\ldots , y_t$ is a d-sequence of $M'$. So we have $c_t-b_{t,1,i_2,\ldots,i_k}\in  IM :_M y_1y_t=IM :_M y_t$. Therefore  $y_t(c_t-b_{t,1,i_2,\ldots,i_k})\in IM$. 
It follows that there exist $d_{t+1},\ldots , d_{s-2}\in M$ such that  
$$y_sa_{1,i_2,\ldots,i_k}=\sum_{j=2}^t y_jb_{j,1,i_2,\ldots,i_k}+\sum_{j=t+1}^{s-2}y_jd_j.$$
 If $t<s-2$, we can repeat  the above process. Finally, after at most $(s-2-t)$-steps, we get
$$y_sa_{1,i_2,\ldots,i_k}=\sum_{j=2}^{s-2} y_jb_{j,1,i_2,\ldots,i_k}$$
as required and the proof of the lemma is complete.
\end{proof}
The following lemma proved in \cite {cm} is needed for the proofs  of the theorem \ref {main} and its consequences.
\begin{lemma} {\rm \cite[Lemma 3.1]{cm}}\label{cmm}
Let $\un x=(x_1,\ldots,x_d)$ be a system of parameters of  $M$. If $\un x$ is a strong d-sequence, then for  $0< i\leq j\leq d$ the length of the Koszul homology module $H_i(x_1^{n_1},\ldots,x_j^{n_j};M)$ is finite and given by
$$\ell(H_i(x_1^{n_1},\ldots,x_j^{n_j};M))=\sum^{j-i}_{t=0}\dbinom{j-t-1}{i-1}\ell(H^0_{\fm}(M/(x_1^{n_1},\ldots,x_t^{n_t})M)).$$
Therefore the function $\ell(H_i(x_1^{n_1},\ldots,x_j^{n_j};M))$  depends only on $n_1,\ldots , n_{j-i}$.
\end{lemma}
 
Now we are able to prove Theorem  \ref {main}.

\noindent
{\it Proof of Theorem \ref {main}}. Since $\un x$ is a dd-sequence,  $\chi_{d-k}(\un x(\un n);M)$ is a polynomial in $\un n$ by Theorem \ref {main3} and Lemma \ref {lm}. Then $\text {deg }\chi_{d-k}(\un x(\un n);M)=p_{d-k}(M)\leq k$ by Remark \ref {r1}. Therefore the theorem is proved if we can show that 
\begin {equation}
\chi_{d-k}(\un x(\un n);M)=\sum_{i=0}^{k}n_1\ldots n_ie(x_1,\ldots,x_i;(0:x_{i+1})_{H_{d-k-1}(x_{i+2},\ldots,x_d;M)}) \tag {\#} 
\end {equation}
for all $k=0,1,\ldots ,d-1$. We will do this by the recursive method.
  Denote by $\un x_{k+1}$ the sequence $x_{k+1},\ldots,x_d$. 
Then by virtue of  Lemma \ref{cmm} and applying  Lemma \ref{ckh} to the sequence $x_1^{n_1},\ldots,x_k^{n_k},x_{k+1},$ $\ldots,x_d$ we get
\[ \begin{aligned}
\chi_{d-k}(\underline x(\underline n);M)=&
\chi_{d-k}(x_1^{n_1},\ldots,x_k^{n_k},\un x_{k+1};M)\\ 
=&n_1\ldots n_ke(x_1,\ldots ,x_k;(0:x_{k+1})_{H_{d-k-1}(\un x_{k+2};M)})\\
+& \sum_{i=0}^{k-1}n_1\ldots n_ie(x_1,\ldots ,x_i;(0:x_{i+1}^{n_{i+1}})_{H_{d-k-1}(x_{i+2}^{n_{i+2}},\ldots ,x_k^{n_k},\un x_{k+1};M)}).
\end{aligned} \]
Put
\[ \begin{aligned}
\phi_1(n_1,\ldots,n_k)=&
\chi_{d-k}(x_1^{n_1},\ldots,x_k^{n_k},\un x_{k+1};M)\\ 
& -n_1\ldots n_ke(x_1,\ldots ,x_k;(0:x_{k+1})_{H_{d-k-1}(\un x_{k+2};M)})\\
=& \sum_{i=0}^{k-1}n_1\ldots n_ie(x_1,\ldots ,x_i;(0:x_{i+1}^{n_{i+1}})_{H_{d-k-1}(x_{i+2}^{n_{i+2}},\ldots ,x_k^{n_k},\un x_{k+1};M)}).
\end{aligned} \]
First, we claim that the function $\phi_1(n_1,\ldots,n_k)$ is independent of $n_k$. In fact,
since the function $\chi_{d-k}(\un x(\un n));M)$ does not depend on the order of $x_1,\ldots,x_d$, using again Lemma \ref{ckh} to the sequence $x_k^{n_k},x_{k-1}^{n_{k-1}},\ldots,x_1^{n_1},x_{k+1},\ldots,x_d$ we obtain
\[ \begin{aligned}
\chi_{d-k}(\un x (\un n));M) &=n_1\ldots n_ke(x_1,\ldots ,x_k;(0:x_{k+1})_{H_{d-k-1}(\un x_{k+2};M)})\\
+& \sum_{i=2}^{k+1}n_i\ldots n_ke(x_i,\ldots ,x_k;(0:x_{i-1}^{n_{i-1}})_{H_{d-k-1}(x_1^{n_1},\ldots ,x_{i-2}^{n_{i-2}},\un x_{k+1};M)}).
\end{aligned}\]
Hence
\[
\begin{aligned}
\phi_1(n_1,\ldots,n_k)=&
\sum_{i=2}^k n_i\ldots n_ke(x_i,\ldots ,x_k;(0:x_{i-1}^{n_{i-1}})_{H_{d-k-1}(x_1^{n_1},\ldots ,x_{i-2}^{n_{i-2}},\un x_{k+1};M)})\\
&+\ell((0:x_k^{n_k})_{H_{d-k-1}(x_1^{n_1},\ldots ,x_{k-1}^{n_{k-1}},\un x_{k+1};M)})).
\end{aligned}
\]
On the other hand, it follows  by Lemma \ref{lm} that \[ \begin{aligned}
(0:x_{i-1}^{n_{i-1}})_{H_{d-k-1}(x_1^{n_1},\ldots ,x_{i-2}^{n_{i-2}},\un x_{k+1};M)})&\subseteq (0:x_{i-1}^{n_{i-1}}x_k)_{H_{d-k-1}(x_1^{n_1},\ldots ,x_{i-2}^{n_{i-2}},\un x_{k+1};M)})\\
&=(0:x_k)_{H_{d-k-1}(x_1^{n_1},\ldots ,x_{i-2}^{n_{i-2}},\un x_{k+1};M)})
\end{aligned} \]
for all $2\leq i\leq k$.
Thus $x_k(0:x_{i-1}^{n_{i-1}})_{H_{d-k-1}(x_1^{n_1},\ldots ,x_{i-2}^{n_{i-2}},\un x_{k+1};M)})=0$. So
 $e(x_i,\ldots ,x_k;(0:x_{i-1}^{n_{i-1}})_{H_{d-k-1}(x_1^{n_1},\ldots ,x_{i-2}^{n_{i-2}},\un x_{k+1};M)})=0$. Therefore, by Lemma \ref{lm},
\[ \begin{aligned}
\phi_1(n_1,\ldots,n_k)&=\ell((0:x_k^{n_k})_{H_{d-k-1}(x_1^{n_1},\ldots ,x_{k-1}^{n_{k-1}},\un x_{k+1};M)})\\
&=\ell((0:x_k)_{H_{d-k-1}(x_1^{n_1},\ldots ,x_{k-1}^{n_{k-1}},\un x_{k+1};M)})
\end{aligned}  \]
 is  independent  of $n_k$. Thus the formula ($\#$) is proved if $k=1$. Assume that $k\geq 2$.
Then
\[ \begin{aligned}
\chi_{d-k}&(\un x(\un n);M)\\    
= & n_1\ldots n_ke(x_1,\ldots ,x_k;(0:x_{k+1})_{H_{d-k-1}(\un x_{k+2};M)})
+\phi_1(n_1,\ldots,n_{k-1},1)\\           
= & n_1\ldots n_ke(x_1,\ldots ,x_k;(0:x_{k+1})_{H_{d-k-1}(\un x_{k+2};M)})\\
& +n_1\ldots n_{k-1}e(x_1,\ldots ,x_{k-1};(0:x_k)_{H_{d-k-1}(\un x_{k+1};M)})\\
& + \sum_{i=0}^{k-2}n_1\ldots n_ie(x_1,\ldots ,x_i;(0:x_{i+1}^{n_{i+1}})_{H_{d-k-1}(x_{i+2}^{n_{i+2}},\ldots ,x_{k-1}^{n_{k-1}},J_k;M)}).
\end{aligned} \]
Next, we set
\[ \begin{aligned}
\phi_2(n_1,&\ldots,n_{k-1})\\
&=\chi_{d-k}(\un x(\un n);M) 
-n_1\ldots n_ke(x_1,\ldots ,x_k;(0:x_{k+1})_{H_{d-k-1}(\un x_{k+2};M)})\\
&\hspace{0.45cm}-n_1\ldots n_{k-1}e(x_1,\ldots ,x_{k-1};(0:x_k)_{H_{d-k-1}(\un x_{k+1};M)})\\
&= \sum_{i=0}^{k-2}n_1\ldots n_ie(x_1,\ldots ,x_i;(0:x_{i+1}^{n_{i+1}})_{H_{d-k-1}(x_{i+2}^{n_{i+2}},\ldots ,x_{k-1}^{n_{k-1}},J_k;M)}).
\end{aligned} \]
With the same method as used above, we can verify that  $\phi_2(n_1,\ldots,n_{k-1})$ is independent of $n_{k-1}$. Continuing this procedure, after k- steps,  we obtain
$$\chi_{d-k}(\un x(\un n);M)=\sum_{i=0}^kn_1\ldots n_ie(x_1,\ldots ,x_i;(0:x_{i+1})_{H_{d-k-1}(x_{i+2},\ldots ,x_d;M)})$$
as required.\qed
\begin {corollary} \label {cl1} Let $\un x=(x_1,\ldots,x_d)$ be a system of parameters of  $M$. If $\un x$ is a  dd-sequence, then 
$\ell (H_k(x_1^{n_1},\ldots,x_d^{n_d};M))$ is a polynomial in $n_1,\ldots , n_{p_k(M)}$ of degree $ p_k(M)$ for all $k=0,1,\ldots , d$.
\end {corollary}
\begin {proof} Since 
$$\chi_{k}(\un x(\un n);M) +\chi_{k+1}(\un x(\un n);M)= \ell (H_k(\un x(\un n);M)),$$
the  statement follows from Theorem \ref {main} and the fact that 
$${\rm deg }( \chi_{k}(\un x(\un n);M))=p_k(M)\geq p_{k+1}(M).$$
\end {proof}
%%%%%%%%%%%%%%
%%%%%%%%%%%%%%%%%%%%%%

\section{Local Cohomology Modules}
\markright{Local cohomology modules}
Recall that a \sop $\un x=(x_1,\ldots , x_d)$  of $M$ is said to be a standard \sop if $\chi_{1}(\un x(\un n);M)$ is a constant for all $n_1,\ldots n_d >0$ and that the module $M$ is a generalized Cohen-Macaulay module if there exists a standard \sop on $M$, i. e. $p_1(M)\leq 0$. It is well-known that $M$ is a generalized Cohen-Macaulay module if and only if $\ell(H^i_{\fm}(M))<\infty$ for all $i=1,\ldots , d-1$.
\begin{lemma}\label{gC}
Let $\un x=(x_1,\ldots , x_d)$ be a standard \sop of  a generalized Cohen-Macaulay module  $M$. Then 
$$\ell(H^i_{\fm}(M))=\sum^i_{j=0}(-1)^{i-j}\dbinom{i}{j}\ell( H^0_{\fm}(M/(x_1,\ldots,x_j)M))$$
for all $ i<d$.
\end{lemma}
\begin{proof}
We proceed by induction on $d$ and $i$. The lemma is trivial  if $d=0,1$ or $i=0$. Assume that $d>1$ and $i>0$. Since $\un x$ is a standard \sop of $M$, $(x_1,\ldots , x_d) H^i_{\fm}(M)=0$ for all $i<d$. Hence we have the following exact sequences
$$0\longrightarrow H^{i-1}_{\fm}(M)\longrightarrow H^{i-1}_{\fm}(M/x_1M)\longrightarrow H^i_{\fm}(M)\longrightarrow 0$$
for all $0<i<d$.
Therefore
$$\ell(H^i_{\fm}(M))=\ell(H^{i-1}_{\fm}(M/x_1M))-\ell(H^{i-1}_{\fm}(M))$$
and the lemma follows easily from the inductive hypotheses.
\end{proof}
A consequence of Theorem \ref{chara} is that if $p_1(M)>0$ and the system of parameters $\un x=(x_1,\ldots , x_d)$ is a dd-sequence of $M$ then $p_1(M/x_1M)=p_1(M)-1$. Therefore, if $k\geq p_1(M)$, $p_1(M/(x_1,\ldots,x_k)M)\leq 0$. So $M/(x_1,\ldots,x_k)M$ is a generalized Cohen-Macaulay module. 
\begin{lemma}\label{dd}
Let $\un x=(x_1,\ldots , x_d)$ be  a system of parameters of $M$  and $k$  a positive integer such that $d-1\geq k\geq p_1(M)$. Assume that $\un x$ is a dd-sequence. Then  $x_{k+1},\ldots,x_d$ is a standard system of parameters of $M/(x_1,\ldots,x_k)M$.
\end{lemma}
\begin{proof} 
Put $M^\prime=M/(x_1,\ldots,x_k)M$. Since $\un x$ is a d-sequence, we have 
\[
\begin{aligned}
\chi_1(x_{k+1},\ldots,x_d;M^\prime)&=\ell(M^\prime/(x_{k+1},\ldots,x_d)M^\prime)-e(x_{k+1},\ldots,x_d;M^\prime)\\
&=\ell(M/\underline xM)-e(\underline x;M)\\
&=\chi_1(\underline x;M).
\end{aligned}
\]
Otherwise,  the function $\chi_1(\un x(\un n);M)$ does not depend on $n_{k+1},\ldots,n_d$ by Theorem \ref{main}. Therefore
\[
\begin{aligned}
\chi_1(x^{n_{k+1}}_{k+1},\ldots,x^{n_d}_d;M^\prime)&=\chi_1(x_1,\ldots,x_k,x^{n_{k+1}}_{k+1},\ldots,x^{n_d}_d;M)\\
&=\chi_1(\un x;M)\\
&=\chi_1(x_{k+1},\ldots,x_d;M^\prime),
\end{aligned}
\]
for all $n_{k+1},\ldots,n_d$. So $x_{k+1},\ldots,x_d$ is a standard \sop of $M^\prime$ as required.
\end{proof}

\noindent
{\it Proof of Theorem \ref {main4}.}
Since the \sop $\un x$ is a dd-sequence, it is not difficult to show  by using Lemma  \ref{cmm} that
$$\ell(H^0_{\fm}(M/(x_1^{n_1},\ldots,x_i^{n_i})M))
=\sum^i_{j=0}(-1)^{i-j}\dbinom{d-j-1}{d-i-1}\ell(H_{d-j}(\underline x(\underline n);M)).$$
Therefore $\ell(H^0_{\fm}(M/(x_1^{n_1},\ldots,x_i^{n_i})M))$ is a polynomial in $\un n$ by Corollary \ref {cl1}.
For  $k\geq p_1(M)$ and  positive integers $n_1,\ldots,n_k$, we set $M_{k(\un n)}=M/(x_1^{n_1},\ldots,x_k^{n_k})M$. Then $M_{k(\un n)}$ is a generalized Cohen-Macaulay module and  $x_{k+1},\ldots,x_d$ is a standard \sop of $M_{k(\un n)}$ by Lemma \ref{dd}. Hence, it follows from  Lemma \ref{gC} that
\[ \begin{aligned}
\ell(H^i_{\fm}(M_{k(\un n)}))&=\sum^i_{j=0}(-1)^{i-j}\dbinom{i}{j}\ell( H^0_{\fm}(M_{k(\un n)}/(x_{k+1},\ldots,x_{k+j})M_{k(\un n)}))\\
&=\sum^i_{j=0}(-1)^{i-j}\dbinom{i}{j}\ell( H^0_{\fm}(M/(x_1^{n_1},\ldots , x_k^{n_k},x_{k+1},\ldots,x_{k+j})M))
\end{aligned} \]
is a polynomial in $n_1,\ldots , n_k$ for all $i=0,1, \ldots , d-k-1$.
\qed
\begin {corollary}
Let $\un x=(x_1,\ldots,x_d)$ be a system of parameters of  $M$. If $\un x$ is a  dd-sequence, then 
$\ell (H_k(x_1^{n_1},\ldots,x_i^{n_i};M))$ is a polynomial in $n_1,\ldots , n_{i-k}$ for all $i=1,\ldots , d$ and $k=1,\ldots , i$.
\end {corollary}
\begin {proof}
We obtain from the proof of Theorem \ref {main4} that $\ell(H^0_{\fm}(M/(x_1^{n_1},\ldots,x_i^{n_i})M))$  is a polynomial in $n_1,\ldots , n_i$. Therefore
 $\ell (H_k(x_1^{n_1},\ldots,x_i^{n_i};M))$ is also a polynomial in   $n_1,\ldots , n_{i-k}$  by Lemma \ref {cmm}.
\end {proof}

\section{Sequentially Cohen-Macaulay modules}
The concept of  sequentially Cohen-Macaulay module was introduced by Stanley \cite {st} for graded modules. Here we define this notion for the local case (see \cite {cn}, \cite {sch1}).

\begin {definition}
 (i) A filtration \dfil  of submodules of $M$ is called the {\it dimension filtration} of $M$ if $M_0=H_\fm^0(M)$ and $M_{i-1}$ is the largest submodule of $M_i$ which has dimension  strictly less than $\dim M_i$ for all $i=1,\ldots ,t.$

\noindent (ii) The module $M$ is said  to be a {\it sequentially Cohen-Macaulay module} if 
each quotient $M_i/M_{i-1}$ is Cohen-Macaulay for all $i=1,\ldots ,t$ .
\end{definition}
\begin {remark}
(i) Because $M$ is Noetherian, the dimension filtration of $M$ always exists and it is unique. Moreover,  let $0=M_0\subset M_1\subset \ldots \subset M_t=M$ be a dimension filtration of $M$ with $\dim M_i=d_i.$ Then we have
$$M_i=\bigcap_{\dim R/\frak p_j > d_{i}}Q_j,$$
for all $i=0, 1, \ldots, t-1,$ where $0=\bigcap_{j=1}^n Q_j$ is a reduced primary decomposition of $0$ in $M$ with $Q_j$ is $\frak p_j-$primary for $j=1,\ldots ,n.$

\noindent
(ii) Let  $0=N_0\subset N_1\subset \ldots \subset N_t=M $ be a filtration of submodules of $M$.  Suppose that 
each quotient $N_i/N_{i-1}$ is Cohen-Macaulay and
 $\dim N_1/N_0<\dim N_2/N_1<\ldots <\dim N_t/N_{t-1}.$ Then it was shown in \cite {cn} (see also \cite {sch1}) that this filtration  is just the dimension filtration of $M$ and therefore $M$ is sequentially Cohen-Macaulay.
\end {remark}

To prove Theorem 1.5 we need some auxiliary lemmata
\begin{lemma}\label{1}
Let $ M_0\subset M_1\subset \cdots \subset M_{t-1}\subset M_t=M$ be the dimension filtration of $M$ with $\dim(M_i)=d_i$ and  $\un x=(x_{1},\ldots,x_d)$ a system of parameters of $M$ which is  a dd-sequence. Then $M_{t-1}=0:_Mx_{d}$. 
\end{lemma}
\begin {proof}
Since $\un x$ is a dd-sequence, $x_d^d \in \frak a(M)$ by Corollary \ref{a(M)}. Therefore $x_d^d \in \Ann_R(0:_Mx)$ for every parameter element $x$ of $M$ by  \cite[Satz 2.4.5]{sch}. By virtue of the Prime Avoidance Theorem we can choose now a parameter element $x$ such that $x\in \Ann_R(M_{t-1})$. Then it follows from the definition of the dimension filtration that
$$ M_{t-1}\subseteq 0:_Mx\subseteq 0:_Mx_d^d= 0:_Mx_d\subseteq M_{t-1}.$$
Thus $M_{t-1}=0:_Mx_{d}$ as required.
\end {proof}
\begin{lemma}\label{2}
Let \dfil be  the dimension filtration of $M$ with $\dim M_i=d_i$. Assume that $M$ is sequentially Cohen-Macaulay and $x\in \mathrm{Ann}(M_{t-1})$ is a parameter element of $M$. Then $M/xM$ is a sequentially \cm module with the dimension filtration
$$\frac{M_0+xM}{xM}\subset \frac{M_1+xM}{xM}\subset \cdots 
\subset \frac{M_s+xM}{xM}\subset \frac{M}{xM},$$
where $s=t-2$ if $d_{t-1}=d-1$ and $s=t-1$ if $d_{t-1}<d-1$.
\end{lemma}
\begin{proof}
Since $x\in \text{Ann}(M_{t-1}), M_{t-1}\subseteq 0:_Mx^n$ for  all $n>0$ and therefore   $M_{t-1}=0:_Mx^n$ by the definition of the dimension filtration. Then it is easy to check that $ M_i\cap xM\subseteq M_{t-1}\cap xM=0:_Mx\cap xM=0$ for all $i\leq t-1$.  So we get
$\frac{M_i+xM}{xM}\cong M_i$. Thus
$\frac{(M_i+xM)/xM}{(M_{i-1}+xM)/xM}\cong  M_i/M_{i-1} $  is  \cm for all $i\leq t-1$.
Therefore $(M_s+xM)/xM$ is sequentially \cm with the dimension filtration
$$\frac{M_0+xM}{xM}\subset \frac{M_1+xM}{xM}\subset \cdots 
\subset \frac{M_s+xM}{xM}$$
by Remark 6.2.
Since $\dim(M_s+xM)/xM=d_s<d$, it remains to prove that the module $\frac{M/xM}{(M_s+xM)/xM}$ is a \cm module.
In fact,  if  $\dim(M_{t-1})<d-1$, i. e. $s=t-1$,     we deduce that 
$$\frac{M/xM}{(M_{t-1}+xM)/xM}\cong \frac {M/M_{t-1}}{x (M/M_{t-1})}$$
is Cohen-Macaulay since $M/M_{t-1}$ is Cohen-Macaulay. 
Assume that $\dim(M_{t-1})=d-1$, i. e. $s=t-2$. We have a short exact sequence,
$$0\longrightarrow \frac{(M_{t-1}+xM)/xM}{(M_{t-2}+xM)/xM}\longrightarrow \frac{M/xM}{(M_{t-2}+xM)/xM}\longrightarrow \frac{M}{M_{t-1}+xM}\longrightarrow 0.$$
where $\frac{(M_{t-1}+xM)/xM}{(M_{t-2}+xM)/xM}\cong M_{t-1}/M_{t-2}$  and $M/(M_{t-1}+xM)\cong  \frac {M_t/M_{t-1}}{x (M_t/M_{t-1})}$ are \cm modules of dimension $d-1$. So $\frac{M/xM}{(M_{t-2}+xM)/xM}$ is a Cohen Macaulay module and
the lemma is proved.
\end {proof}
\noindent 
{\it Proof of Theorem \ref {5}}

\noindent
$(a):$ The implication $  (i) \Rightarrow (ii)$  is trivial.

\noindent
$(ii) \Rightarrow (iii)$ By the assumption $x_1,\ldots,x_{d_i}$ is a \sop of $M_i$, so is of $\overline M_i = M_i/M_{i-1}$.  Since $\overline M_i$ is  Cohen-Macaulay,   we get by \cite [Theorem 5.2, a)] {sch1} that
$$\ell(M/\un xM)=\sum_{i=0}^t\ell(\overline M_i/\un x\overline M_i)=\sum_{i=0}^te(x_1,\ldots,x_{d_i};\overline M_i)
=\sum_{i=0}^te(x_1,\ldots,x_{d_i};M_i).$$
So $$\chi_1(\un x;M)=\sum_{i=0}^{t-1}e(x_1,\ldots,x_{d_i};M_i)$$
and (iii) follows since the \sop  $\un x(\un n)=(x_1^{n_1},\ldots,x_d^{n_d})$ is also satisfied the assumption of (ii)\\
$(iii) \Rightarrow (iv)$ is straightforward by Theorem \ref{chara}.\\
$(iv) \Rightarrow (i)$ 
We proceed by induction on $d=\dim M$. There is nothing to do for the case $d=1$  since $\un x=(x_1)$ is a strong d-sequence. Assume that $d>1$. Then $M_{t-1}=0:_M x_d$  by  Lemma 6.3. Therefore $x_d\in \Ann M_{t-1}$ and $M_{t-1}\cap x_dM=0$.
 Put $M^\prime=M/x_dM$. We obtain by Lemma 6.4 that $M^\prime=M/x_dM$ is sequentially \cm with the dimension filtration
$$\frac{M_0+x_dM}{x_dM}\subset \frac{M_1+x_dM}{x_dM}\subset \cdots \subset \frac{M_{s}+x_dM}{x_dM}\subset \frac{M}{x_dM}=M^\prime, \ (*)$$
where $s=t-2$ if $d_{t-1}=d-1$ and $s=t-1$ if $d_{t-1}<d-1$. Let $d_{t-1}<d-1$, i. e. $s=t-1$ in the filtration (*).
Since  $(x_1,\ldots, x_{d-1})$ is a dd-sequence of $M^\prime$, it follows
 from the inductive hypotheses for all $i\leq t-1$ that
$$(M_i+x_dM)\cap (x_{d_i+1},\ldots, x_{d-1},x_d)M=x_dM.$$
 Thus  
$$M_i\cap(x_{d_i+1},\ldots, x_d)M\subseteq M_i\cap x_dM\subseteq M_{t-1}\cap x_dM=0$$
and the implication follows.  Since $M_{t-1}\cap x_dM=0$, the statement is clear for the case $d_{t-1}=d-1$, i. e. $s=t-2$ in the filtration (*).\\
 \medskip
$(b)$: From the short exact sequence
$0\longrightarrow M_{t-1}\longrightarrow M\longrightarrow  M/M_{t-1}\longrightarrow 0,$
we get a long exact sequence of Koszul homology modules
\begin{multline*}
\cdots \longrightarrow H_{k+1}(\un x; M/M_{t-1}) \longrightarrow H_k(\un x; M_{t-1}) \longrightarrow \\ 
H_k(\un x; M) \longrightarrow  H_k(\un x; M/M_{t-1})\longrightarrow \cdots.
\end{multline*} 
 Since $M/M_{t-1}$ is Cohen-Macaulay,  $H_k(\un x;  M/M_{t-1})=0$ for all $k>0$. Therefore  $H_k(\un x; M_t)\cong H_k(\un x;M_{t-1}).$ So $\chi_k(\un x;M)=\chi_k(\un x; M_{t-1})$. 
For a positive integer $n\leq d$ we set $\un x_n=(x_1,\ldots, x_n)$. Since $x_nM_{t-1}=0$ for all $n>d_{t-1}$ by the hypothesis, $x_nH_i(\un x_{n-1}; M_{t-1})=0$ for  all $i$. Thus, from the long exact sequence
$$\cdots \longrightarrow H_{i+1}(\un x_n;M_{t-1})\longrightarrow H_i(\un x_{n-1};M_{t-1})\stackrel{*x_n}{\longrightarrow} H_i(\un x_{n-1};M_{t-1})\longrightarrow \cdots$$
we get the short exact sequences,
$$0\longrightarrow H_{i+1}(\un x_{n-1}; M_{t-1}) \longrightarrow  H_{i+1}(\un x_n; M_{t-1})\longrightarrow H_i(\un x_{n-1}; M_{t-1}) \longrightarrow 0.$$
Therefore we can show by recursive method on $n=d,\ldots , d_{t-1}+1$ that  
$$\ell(H_i(\un x; M_{t-1})) =\sum^i_{j=0}\binom{d_t-d_{t-1}}{j}\ell(H_{i-j}(\un x_{d_{t-1}}; M_{t-1}))$$ for all $i\geq 1$.
Then for  $k>0$ we have
\[ \begin{aligned}
\chi_k(\un x;M)&=\sum^{d_t}_{i=k}(-1)^{i-k}\ell(H_i(\un x;M_{t-1}))\\
&=\sum^{d_t}_{i=k}(-1)^{i-k}\sum^i_{j=0}\binom{d_t-d_{t-1}}{j}\ell(H_{i-j}(\un x_{d_{t-1}}; M_{t-1}))\\
&=\sum^{k-1}_{j=0}\binom{d_t-d_{t-1}}{j}\sum^{d_t}_{i=k}(-1)^{i-k}\ell(H_{i-j}(\un x_{d_{t-1}}; M_{t-1}))\\
&\hspace{0.5cm}+\sum^{d_t}_{j=k}\binom{d_t-d_{t-1}}{j}\sum^{d_t}_{i=j}(-1)^{i-k}\ell(H_{i-j}(\un x_{d_{t-1}}; M_{t-1}))\\
&=\sum^{k-1}_{j=0}\binom{d_t-d_{t-1}}{j}\chi_{k-j}(\un x_{d_{t-1}};M_{t-1})\\
&\hspace{0.5cm}+
(-1)^k\big(\sum^{d_t}_{j=k}(-1)^j\binom{d_t-d_{t-1}}{j}\big)e(\un x_{d_{t-1}};M_{t-1}).
\end{aligned} \]
Now, we can obtain the formular of  the statement (b) by induction on $t$. 
Hence  the proof of Theorem 1.5 is complete. \qed

Note that a \sop $\un x=(x_1,\ldots , x_d)$, which is satisfied the condition (ii) of Theorem 1.5, is just called by P. Schenzel  (see \cite [Definition 2.5] {sch1}) a distinguished system of parameters. Since there always  exists a distinguished \sop on any finitely generated module, it follows immediately from Theorem 1.5, (a) the existence of system of parameters, which is a dd-sequence, on a sequentially \cm module.
\begin{corollary}
Any sequentially \cm module admits a system of parameters which is a dd-sequence.
\end {corollary}

 An other immediate consequence of Theorem 1.5, (b) is to present a relation between the invariants $p_k(M), (0\leq k\leq d)$ and the dimension $d_i=\dim M_i$ as follows.

\begin{corollary}
Suppose that $M$ is a sequentially \cm module. Then for all $i=0,\ldots ,d$
 $$p_k(M)=d_i=\dim(R/\fa_{d_i}(M))$$
for  $d-d_{i+1}<k\leq d- d_{i}$, 
where $\fa_k(M)=\mathrm{Ann}(H^k_\fm(M))$.
\end{corollary}

\vspace{1cm}

%%%%%%%%%%%%%%%%%%%%%%%

\end {document}